\documentclass[11pt]{article}
\usepackage{amsthm, amsmath, amssymb, amsfonts, url, booktabs, tikz, setspace, fancyhdr, enumerate}
\usepackage[margin = 1in]{geometry}

\usepackage[czech,english]{babel} 

\usepackage[hidelinks]{hyperref}
\usepackage{tikz}
\usepackage{comment}
\usepackage{xcolor}
\usepackage{cite}

\usepackage[textsize=scriptsize,colorinlistoftodos]{todonotes}

\newtheorem{theorem}{Theorem}[section]

\newtheorem{lemma}[theorem]{Lemma}
\newtheorem{corollary}[theorem]{Corollary}

\newtheorem{claim}[theorem]{Claim}

\theoremstyle{definition}

\theoremstyle{remark}

\newcommand{\Hom}{\mathrm{Hom}}

\newcommand{\HH}{\mathcal{H}}

\newcommand{\C}{\mathcal{C}}
\newcommand{\D}{\mathcal{D}}

\newcommand{\F}{\mathcal{F}}
\newcommand{\Oh}{\mathcal{O}}

\newcommand{\U}{\mathcal{U}}

\newcommand{\ex}{\mathrm{ex}}

\title{Rainbow cycles in properly edge-colored graphs}
\author{Jaehoon Kim\thanks{Department of Mathematical Sciences, KAIST, South Korea. Email: \texttt{jaehoon.kim}@\texttt{kaist.ac.kr}. Supported by the POSCO Science Fellowship of POSCO TJ Park Foundation}
\and Joonkyung Lee\thanks{Department of Mathematics, 
		Hanyang University, 
		222 Wangsimni-ro, Seongdong-gu, 
		Seoul 04763, South Korea. 
		Email: \texttt{joonkyunglee}@\texttt{hanyang.ac.kr}.}
\and Hong Liu
	\thanks{Extremal Combinatorics and Probability Group (ECOPRO), Institute for Basic Science (IBS), Daejeon, South Korea. Email: {\tt hongliu@ibs.re.kr}. Supported by IBS-R029-C4.}	
\and Tuan Tran\thanks{School of Mathematical Sciences, University of Science and Technology of China, China. E-mail: {\tt trantuan@ustc.edu.cn}. Supported by the Outstanding Young Talents Program (Overseas) of the National Natural Science Foundation of China.}
 }
\date{}

\begin{document}

\maketitle

\begin{abstract}
We prove that every properly edge-colored $n$-vertex graph with average degree at least $100(\log n)^2$ contains a rainbow cycle, improving upon 
$(\log n)^{2+o(1)}$ bound due to Tomon. 
We also prove that every properly colored $n$-vertex graph with at least $10^5 k^2 n^{1+1/k}$ edges contains a rainbow $2k$-cycle, which improves the previous bound $2^{ck^2}n^{1+1/k}$ obtained by Janzer.

Our method using homomorphism inequalities and a lopsided regularization lemma also provides a simple way to prove the Erd\H{o}s--Simonovits supersaturation theorem for even cycles, which may be of independent interest.
\end{abstract}

\section{Introduction}
In an edge-colored graph, a copy of a graph $H$ is called \emph{rainbow} if every edge in the copy receives a unique color.
If we forbid rainbow copies of some graphs in properly edge-colored graphs~$G$, what is the maximum number of edges in $G$?
This extremal question was first investigated by
Keevash, Mubayi, Sudakov, and Verstra\"ete~\cite{KMSV07},
where they defined the \emph{rainbow Tur\'an number} $\ex^*(n,\HH)$ for a family of graphs $\HH$.
Formally, $\ex^*(n,\HH)$ denotes the maximum number of edges in a properly edge-colored $n$-vertex graph with no rainbow copies of any $H\in \HH$.
If $\HH$ consists of a single graph $H$, then we simply write $\ex^*(n,H)$ instead of $\ex^*(n,\{H\})$.

To quote~\cite{KMSV07}, there are two questions that are the most important among the several ones raised therein. The first one is to determine $\ex^*(n,\C)$, where $\C$ is the class of all cycles. It is shown that $\ex^*(n,\C)=\Omega(n\log n)$ in~\cite{KMSV07} and Das, Lee, and Sudakov~\cite{DLS13} obtained an upper bound $O(ne^{(\log n)^{\frac{1}{2}+o(1)}})$.
There have been some recent improvements upon the upper bound~\cite{J20,T22}  
and the current best one is $O(n(\log n)^{2+o(1)})$ appeared in \cite{T22}. 
We improve this bound to $O(n\log^2 n)$.
\begin{theorem}\label{thm:cycle}
A properly edge-colored $n$-vertex graph $G$ 
with 
at least $32n\log^2 (5n)$ edges always contains a rainbow cycle.
\end{theorem}

Theorem~\ref{thm:cycle} has the following corollary on cycles in linear 3-uniform hypergraphs.
\begin{corollary}\label{cor:3-graph}
If an $n$-vertex linear $3$-uniform hypergraph $H$ has at least $288 n \log^2(5n)$ edges, it contains a loose cycle.
\end{corollary}

The second question in~\cite{KMSV07} concerns with $\ex^*(n,C_{2k})$, where $C_{2k}$ is the even cycle of length $2k$. In~\cite{KMSV07}, a general lower bound $\ex^*(n,C_{2k})=\Omega(n^{1+1/k})$ is obtained, whereas the matching upper bounds were only verified for $k=2,3$. This upper bound was subsequently improved by Das, Lee, and Sudakov~\cite{DLS13} to $O(n^{1+(1+o_k(1))\log k/k})$ and by Janzer~\cite{J20} to $O(n^{1+1/k})$.
While Janzer's bound matches the lower bound given in~\cite{KMSV07}, the implicit constant is exponential in $k$. We improve it to a polynomial one as follows.
\begin{theorem}\label{thm:even}
A properly edge-colored $n$-vertex graph $G$ with at least $10^5 k^3 n^{1+1/k}$ edges always contains a rainbow $2k$-cycle.
\end{theorem}
Whereas the 
recent improvement in~\cite{T22} took a different approach to that of Janzer~\cite{J20}, our proof is in spirit closer to that of~\cite{J20} using homomorphism counts and improves it in two ways.
First, we use more efficient lopsided regularization lemma than the Jiang--Seiver lemma~\cite{JS12} used by Janzer.
Second, the main proof after regularization is conceptually simpler and more intuitive in the sense that it only relies on repeated applications of the Cauchy--Schwarz inequality.

We stress that the repeated Cauchy--Schwarz method may be of independent interest.
There are various problems in extremal combinatorics, from the classical Mantel's theorem and the K\H{o}vari--S\'os--Tur\'an theorem to the recent developments on Sidorenko's conjecture, where the Cauchy--Schwarz inequality and its variants have been extremely useful;
however, the convexity inequalities have seen less success in determining extremal numbers of bipartite graphs other than complete bipartite graphs. 

To the best of our knowledge, it has been unknown whether even the Bondy--Simonovits theorem \cite{BS74}, a weaker statement than Theorem~\ref{thm:even}, has a proof that only uses Cauchy--Schwarz inequality.
Our proof method answers this natural question by giving a simple proof of the Bondy--Simonovits theorem and moreover, obtains the Erd\H{o}s--Simonovits supersaturation \cite{ESi84} as follows.
\begin{corollary}\label{cor:supsat}
  Any $n$-vertex graph $G$ with average degree $d\ge 2\cdot 10^5 k^3 n^{1/k}$ contains at least $\frac{1}{2}(2^{12}k)^{-k}d^{2k}$ copies of $2k$-cycle.
\end{corollary}

\section{The key homomorphism inequality}
In what follows, a \emph{coloring} always means an edge coloring and likewise, a \emph{colored graph} is an edge-colored graph.
Let $\Hom(H,G)$ be the set of homomorphisms from $H$ to $G$. In particular, we fully label vertices of $H$.

Each homomorphism in $\Hom(C_{2k},G)$ can be seen as 
a closed walk $v_0v_1\dots v_{2k}$ of length~$2k$ (or closed $2k$-walk for short) in $G$ with $v_0=v_{2k}$. 
Let $\phi$ be a proper coloring of a graph~$G$.
A closed $2k$-walk is \emph{degenerate} if, for distinct $i$ and $j$ with $\{i,j\}\neq \{0,2k\}$,  $v_i=v_j$ or $\phi(v_iv_{i+1}) = \phi(v_jv_{j+1})$, where the index addition is taken modulo $2k$.
That is, the walk revisits a vertex in the middle or repeats a color.
The former case is said to be \emph{vertex-degenerate at $(i,j)$} and the latter is called \emph{color-degenerate at $(i,j)$}.
If a $2k$-walk is vertex-degenerate (resp. color-degenerate) at $(i,j)$, then it is of \emph{type} $|i-j|-1$ (resp.~$|i-j|$). 

One degenerate walk may have multiple types, so the types of degeneracy are not disjoint in general.
As $\phi$ is a proper coloring, 
a closed walk is vertex-degenerate of type $1$ if and only if it is color-degenerate of type $1$. 
Note that a non-degenerate $2k$-walk is a rainbow cycle of length $2k$.
Let $\D(C_{2k},G)$ be the set of all degenerate closed walks of length $2k$ and let $S_k$ be the star with $k$ leaves. 
If $G$ has no rainbow $2k$-cycle, then $\Hom(C_{2k},G)= \D(C_{2k},G)$. 
Our strategy is to bound $|\D(C_{2k},G)|$ from above to obtain an upper bound for $|\Hom(C_{2k},G)|$.
 
\begin{lemma}\label{lem: degenerate}
Suppose that $G$ has no rainbow $2k$-cycle. Then
$$|\Hom(C_{2k},G)| \leq (2k)^{2k} |\Hom(S_k,G)|.$$
\end{lemma}
\begin{proof}
Let $\U_s$ and $\F_s$ be the set of closed $2k$-walks that are vertex-degenerate at $(0,s+1)$ and color-degenerate at $(0,s)$, respectively, i.e., they consist of those walks of type $s$. Note that $\U_1=\F_1$.
\begin{claim}\label{cl: 1}
For $1\leq s\leq 2k-3$ and $1\leq t \leq 2k-1$,
\begin{align*}
    |\U_s|^2\leq |\U_1| \cdot |\U_{2s-1}| \enspace \text{and} \enspace |\F_t|^2 \leq |\F_1|\cdot |\F_{2t-1}|.
\end{align*}
In particular, $|\U_s|\leq|\U_1|$ and $|\F_t|\leq |\F_1|$.
\end{claim}
\begin{proof}[Proof of the claim]
For $s=1$, the inequality becomes a trivial equality and, for $1\leq s\leq k-1$, $|\U_{s}|=|\U_{2k-2-s}|$ by symmetry.
We may thus assume $1<s\leq k-1$.
We analyze $\U_s$ by counting the closed walks therein in the following way:
\begin{itemize}
    \item fix vertices $x=v_0=v_{s+1}$, $y=v_{s+k}$, and $z=v_s$ with $xz\in E(G)$;
    \item choose a walk $v_{s+1}v_{s+2}\dots v_{s+k}$ of length $k-1$ from $x$ to $y$;
    \item choose a walk $v_{s+k}v_{s+k+1}\dots v_{2k}$ of length $k-s$ from $y$ to $x$; and
    \item choose a walk $v_0v_1\dots v_{s}$ of length $s$ from $x$ to $z$.
\end{itemize}
Let $w_{\ell}(u,v)$ be the number of walks of length $\ell$ from $u$ to $v$
and let $g(u,v)$ be the edge indicator function of $G$, i.e., $g(u,v)=1$ if $uv\in E(G)$ and $0$ otherwise.
Then 
\begin{align*}
    |\U_s| = \sum_{x,y,z} w_{k-1}(x,y) w_{k-s}(y,x)w_{s}(x,z)g(x,z).
\end{align*}
The Cauchy--Schwarz inequality now gives the following bound:
\begin{align*}
    |\U_s|^2
    &\leq \left(\sum_{x,y,z} w_{k-1}(x,y)^2g(x,z)\right) \left(\sum_{x,y,z} w_{k-s}(y,x)^2w_{s}(x,z)^2 g(x,z)\right) \\
    &\leq \left(\sum_{x,y,z} w_{k-1}(x,y)^2g(x,z)\right) \left(\sum_{x,y,z} w_{k-s}(y,x)^2w_{s}(x,z)^2\right),
\end{align*}
where the second inequality follows from $g(x,y)\leq 1$. 
The sum $\sum_{x,y,z} w_{k-1}(x,y)^2g(x,z)$ counts the number of closed $(2k-2)$-walks plus a pendant edge that corresponds to~$xz$, which are exactly those $2k$-walks in $\U_1$. Thus, $\sum_{x,y,z} w_{k-1}(x,y)^2g(x,z)=|\U_1|$. 
Analogously, the term $w_{k-s}(y,x)^2$ and $w_{s}(x,z)^2$ count 
the number of closed walks of length $2(k-s)$ and $2s$ when summed over the choices of $y$ and $z$, respectively, both starting at $x$.
Thus, summing  $w_{k-s}(y,x)^2w_{s}(x,z)^2$ over $x,y,z\in V(G)$ counts the number of closed $2k$-walks that are vertex-degenerate at $(0,2s)$, which is exactly $|\U_{2s-1}|$.
Therefore, $|\U_s|^2\leq |\U_1|\cdot|\U_{2s-1}|$. 
Let $\U_j$ be the largest set among $\U_1,\cdots,\U_{k-1}$. 
Then $|\U_j|\leq |\U_1|^{1/2}|\U_{2j-1}|^{1/2}\leq |\U_1|^{1/2}|\U_j|^{1/2}$,
so $|\U_j|=|\U_1|$. This proves $|\U_s|\leq|\U_1|$ for all $s=1,2,\cdots,2k-3$.

\medskip

The second inequality for $\F_t$ can also be obtained by using essentially the same technique. 
Again by symmetry, $|\F_t|= |\F_{2k-t}|$ for $1\leq t\leq k$, so we may assume $k\leq t< 2k-1$.
As each closed $2k$-walk $v_0\dots v_{2k}$ in $\F_t$ satisfies $\phi(v_0v_{1}) = \phi(v_{t}v_{t+1})$, we count $|\F_t|$ as follows: 
\begin{itemize}
    \item fix a color $c$ that repeats at $v_0v_1$ and $v_{t}v_{t+1}$;
    \item fix vertices $x=v_0=v_{2k}$ and $y=v_{k}$;
    \item choose a walk $v_{0}v_1\dots v_{k}$ of length $k$ from $x$ to $y$ where $\phi(v_0v_{1})=c$; and
    \item choose a walk $v_{k}v_{k+1}\dots v_{2k}$ of length $k$ from $y$ to $x$ where $\phi(v_{t}v_{t+1})=c$. 
\end{itemize}
Let $\tilde{w}_{k,\ell}(u,v,c)$ be the number of $k$-walks from $u$ to $v$ such that 
the $\ell$-th edge has color $c$.
The Cauchy--Schwarz inequality then gives
\begin{align*}
    |\F_t|^2 &= \left(\sum_{x,y,c} \tilde{w}_{k,1}(x,y,c) \tilde{w}_{k,t-k+1}(y,x,c)\right)^2 \\
    &\leq\left(\sum_{x,y,c} \tilde{w}_{k,1}(x,y,c)^2\right)\left(\sum_{x,y,c} \tilde{w}_{k,t-k+1}(y,x,c)^2\right).
\end{align*}
The sum $\sum_{x,y,c} \tilde{w}_{k,1}(x,y,c)^2$ counts the number of closed $2k$-walks from $x$ that repeat $c$ at the first and the last edges, which is exactly $|\F_1|$ by rotational symmetry.
Similarly, the second sum corresponds to $|\F_{2t-2k+1}|=|\F_{2t-1}|$.
Therefore, we obtain the inequality  $|\F_t|^2 \leq |\F_1|\cdot |\F_{2t-1}|$.
Finally, $|\F_t|\leq |\F_1|$ follows from the same argument used for showing $|\U_s|\leq |\U_1|$.
\end{proof}

\begin{claim}\label{cl: 2}
 $|\Hom(C_{2k},G)| \leq 4k^2 |\U_1|$.
\end{claim}
\begin{proof}[Proof of the claim]
Recall that both $\U_s$ and $\F_s$ specify the labels of the vertices where degeneracy occurs. 
By rotational symmetry, the number of closed $2k$-walks that are vertex-degenerate at $(i,i+s+1)$ is exactly $|\U_s|$ for each $i=0,1,\dots, 2k-1$. Likewise, the number of closed $2k$-walks that are color-degenerate at $(i,i+s)$ is $|\F_s|$ for each $i=1,2,\dots,2k$. Thus, the number of degenerate $C_{2k}$-homomorphisms of type $s$ is at most $2k(|\U_s|+|\F_s|)$.
Taking the union bound over all possible types, we get
\begin{align}
    |\D(C_{2k},G)| \leq 2k\left(\sum_{s=1}^{k-1}|\U_s| + \sum_{s=1}^k |\F_s|\right).
\end{align}
Together with Claim~\ref{cl: 1} this yields $|\D(C_{2k},G)| \leq 4k^2|\U_1|$ as desired.
\end{proof}

Let $\Oh_s$ be the set of all closed $2k$-walks $v_0\dots v_{2k}$ with $v_0= v_{2k}$
such that $v_0 = v_2 = \dots = v_{2s}$.
In particular, $\Oh_1 = \U_1 = \F_1$. We also use $\Oh_0 = \Hom(C_{2k},G)$ for consistency.
\begin{claim}\label{cl: 3}
The sequence $|\Oh_s|$, $0\leq s\leq k$, is log-convex, i.e.,  
for each $s=1,\dots, k-1$,
\begin{align*}
     |\Oh_{s}|^2\leq |\Oh_{s-1}|\cdot|\Oh_{s+1}|.
\end{align*}
\end{claim}
\begin{proof}[Proof of the claim]
A \emph{star-walk} of length $\ell$ is a walk $u_0u_1\dots u_{\ell}$ of length $\ell$ such that every even-indexed vertex is the same one, i.e., $u_0=u_2=\dots=u_{2t}$ where $t=\lfloor\tfrac{\ell}{2}\rfloor$.
For $1\leq s\leq k-1$, the walks in $\Oh_s$ can be counted as follows:
\begin{itemize}
    \item fix vertices $x=v_0=v_2=\dots=v_{2s}$, $y=v_{s+1}$, $z=v_{k+s+1}$, where $y$ is either $x$ or a neighbor of $x$ depending on the parity of $s$;
    \item \label{it:walk1} choose a walk $v_{k+s+1}v_{k+s+2}\dots v_{2k}$ of length $k-s-1$ from $z$ to $x$;
    \item \label{it:star1} choose a star-walk $v_0v_1\dots v_{s+1}$ of length $s+1$ from $x$ to $y$;
    \item \label{it:walk2} choose a walk $v_{k+s+1}v_{k+s}\dots v_{2s}$ of length $k-s+1$ from $z$ to $x$; and
    \item \label{it:star2} choose a star-walk $v_{2s}v_{2s-1}\dots v_{s+1}$ of length $s-1$ from $x$ to $y$;
\end{itemize}
As in the proof of Claim~\ref{cl: 1}, $w_{\ell}(u,v)$ denotes the number of $\ell$-walks from $u$ to $v$.
Let $\sigma_{\ell}(u,v)$ be the number of star-walks of length $\ell$ from $u$ to $v$. Note that, unlike $w_\ell(\cdot,\cdot)$, $\sigma_\ell(\cdot,\cdot)$ may not be a symmetric function in general.
Now we can compute $|\Oh_s|$ as
\begin{align*}
    |\Oh_s| = \sum_{x,y,z} w_{k-s-1}(z,x)\sigma_{s+1}(x,y)w_{k-s+1}(z,x)\sigma_{s-1}(x,y).
\end{align*}
The Cauchy--Schwarz inequality then yields 
\begin{align*}
    |\Oh_s|^2 \leq\left(\sum_{x,y,z} w_{k-s-1}(z,x)^2\sigma_{s+1}(x,y)^2\right)\left(\sum_{x,y,z} w_{k-s+1}(z,x)^2\sigma_{s-1}(x,y)^2\right).
\end{align*}
In the first sum, $\sigma_{s+1}(x,y)^2$ counts either a closed star-walk of length $2(s+1)$ (if $s$ is odd and $x=y$) from $x$ or a star-walk of length $2s$ together with an edge $xy$ (if $s$ is even and $xy\in E(G)$). By considering $xy$ and $yx$ as the $(s+1)$-th and $(s+2)$-th edge of the star-walk, summing the latter case over all choices of $y$ counts the closed star-walks of length $2(s+1)$.
Hence the first sum counts the number of walks in $\Oh_{s+1}$ by summing the number of ways to augment each closed walk of length $2(k-s-1)$ and to a closed star-walk of length $2(s+1)$ at $x$.
By the same reason with replacing $s+1$ by $s-1$, the second sum counts the number of walks in $\Oh_{s-1}$, which proves the claim.
\end{proof}
We are now ready to finish the proof of the lemma.
By Claims~\ref{cl: 2} and~\ref{cl: 3}, 
\begin{align*}
    4k^2\geq \frac{|\Hom(C_{2k},G)|}{|\U_1|} = \frac{|\Oh_0|}{|\Oh_1|} \geq \frac{|\Oh_1|}{|\Oh_2|} \geq \dots \geq \frac{|\Oh_{k-1}|}{|\Oh_k|}.
\end{align*}
Therefore, 
\begin{align*}
   \frac{|\Oh_0|}{|\Oh_k|}  = \frac{|\Oh_0|}{|\Oh_1|} \cdot  \frac{|\Oh_1|}{|\Oh_2|}  \dots \frac{|\Oh_{k-1}|}{|\Oh_{k}|}  \leq (4k^2)^k.
\end{align*}
Each closed $2k$-walk in $\Oh_k$ is of the form $vu_1vu_2\dots vu_kv$, which naturally corresponds to a homomorphism from $S_k$ to $G$ that maps the central vertex to $v$ and the $i$-th leaf to $u_i$.
Therefore, $|\Hom(C_{2k},G)|=|\Oh_0|\leq 4^k k^{2k}|\Oh_k|=|\Hom(S_{k},G)|$, which concludes the proof.
\end{proof}

\section{Regularization}
Suppose that the $n$-vertex graph $G$ is ``close" to being regular, e.g., $|\Hom(S_k,G)|\leq 10n^{k+1}p^k$, where $p=2e(G)/n^2$ is the edge density of $G$.
Here the constant $10$ is arbitrarily chosen to illustrate.
If~$G$ contains no rainbow $2k$-cycles, then Lemma~\ref{lem: degenerate} gives
\begin{align*}
    n^{2k}p^{2k}\leq |\Hom(C_{2k},G)|\leq  (2k)^{2k}|\Hom(S_k,G)|\leq 10\cdot (2k)^{2k} n^{k+1}p^{k},
\end{align*}
where the first inequality follows from the fact that even cycles satisfy Sidorenko's conjecture. As a corollary, $e(G)\leq 2\cdot10^{1/k}k^2n^{1+1/k}$, which is stronger than Theorem~\ref{thm:even}.


However, this ideal assumption is not guaranteed in general. Instead, the following two lemmas will enable us to ``regularize" the graph $G$.

\begin{lemma}\label{lem: regularization 2}
Let $G$ be a properly edge-colored graph with minimum degree $\delta(G)\geq 1$. Then there exists a properly edge-colored graph $G'$ with the following properties:
\begin{itemize}
    \item[\rm (1)] $G'$ satisfies $|V(G)|\leq |V(G')|\leq 4e(G)/\delta(G)$ and every vertex of $G'$ has degree between $\delta(G)/2$ and $\delta(G)$;
    \item[\rm(2)] There is a color-preserving homomorphism $\psi$ from $G'$ to $G$. In particular, if $G'$ contains a rainbow cycle, then so does $G$.
\end{itemize}
\end{lemma}
\begin{proof}
Let $\delta:=\delta(G)$ for brevity. 
We construct $G'$ by iterating the following process. Fix an ordering the vertices of $G$. At each step, take a vertex $v\in V(G)$ according to the ordering and let $s:=\lceil d_G(v)/\delta\rceil$. 
We then split $v$ into new vertices $v_1,\dots, v_s$ so that the neighbor sets $N(v_1),\dots, N(v_s)$ form a partition of the neighbor set of $v$ in $G$ and $\delta/2 \leq |N(v_i)|\leq \delta$ for every $1\le i \leq s$. 
This is possible since $d_G(v)\ge \delta$ and we can make the sizes of $N(v_1),N(v_2),\cdots, N(v_s)$ as equal as possible. We color the edges in such a way that each edge $uv_i$ for $u\in N_G(v)$ inherits the same color as $uv$.

Let $\psi$ be a map from $G'$ to $G$ so that each vertex maps to the original vertex of $G$ before splitting. Then $\psi$ is a color-preserving homomorphism and every vertex of $G'$ has degree between $\delta/2$ and $\delta$. Since the number of edges is preserved throughout the whole process,
\[
|V(G')|\cdot \delta/2 \leq 2e(G')=2e(G),
\]
which implies $|V(G')|\leq 4e(G)/\delta$. Indeed, the color-preserving homomorphism $\psi$ maps a rainbow cycle in~$G'$ to a rainbow cycle in $G$. This proves the ``in particular" part.
\end{proof}

\begin{lemma}\label{lem: regularization 1}
Let $k\ge 2$ and let $G$ be an $n$-vertex bipartite graph with 
average degree 
$d>0$.
Suppose that $G$ contains no proper subgraph with larger average degree.
Then $G$ contains a subgraph $G'$ with bipartition $(A,B)$ satisfying the following for some $i\in \mathbb{N}$:
\begin{itemize}
    \item[\rm(1)] $|A| \ge \frac{1}{k}2^{-\frac{ki}{k-1}}n$, $|B|\ge \frac{n}{64}$;
    \item[\rm(2)] $d_{G'}(a) \in [2^{i-6} d, 2^{i-5} d ]$ for all $a\in A$;
    \item[\rm(3)] $d_{G'}(b) \leq 4d$ for all $b\in B$.
\end{itemize}
\end{lemma}
\begin{proof}
Denote by $(X,Y)$ a bipartition of $G$. 
Let $X_0$ and $Y_0$ be the set of vertices in $X$ and $Y$, respectively, of degree at least $4d$ and let $X_1:=X\setminus X_0$ and $Y_1:=Y\setminus Y_0$.
Then $|X_0|, |Y_0|\leq e(G)/(4d)=n/8$.


Since there is no proper subgraph of $G$ with average degree larger than $d$, 
\[
e(G[X_0,Y_0]) \le \frac12(|X_0|+|Y_0|)d \le \frac{nd}{8}=\frac{1}{4}e(G).
\]
Hence, one of $G[X_0,Y_1]$, $G[X_1, Y_0]$ and $G[X_1,Y_1]$ has at least $e(G)/4$ edges. By symmetry, we can assume that for some $X'\in \{X_0,X_1\}$, the graph $G[X',Y_1]$ has at least $e(G)/4$ edges. 

As 
$G[X',Y_1]$ has at least $e(G)/4$ edges and average degree at least $d/4$, we can delete some vertices of 
$G[X',Y_1]$
to obtain a graph $G_1=G[X^*,Y^*]$ with minimum degree at least $d/16$ and $e(G_1)\ge e(G)/8$. 
Since the vertices in $Y^*$ have degree at most $4d$,
\[
|Y^*|\ge \frac{e(G_1)}{4d} \ge \frac{e(G)}{32d}=\frac{n}{64}.
\]


We partition 
$X^*$ into the following sets 
\[
Z_i = \left\{ v\in X^*\colon d_{G_1}(v) \in [2^{i-6}d, 2^{i-5}d) \right\}, \quad i\in \mathbb{N}.
\]
If there exists $i$ such that $|Z_i| \geq \frac{1}{k}2^{-\frac{ki}{k-1}}n$, then take $A$ to be 
$Z_i$ and $B=Y^*$, and we are done.
If not, then we have
$$\frac{e(G)}{8} \leq e(G_1) \le \sum_{i=1}^{\infty} |Z_i|\cdot 2^{i-5}d  < \frac{dn}{32k}\cdot\sum_{i=1}^{\infty} 2^{-\frac{i}{k-1}}  < \frac{d n}{32 k} \cdot \frac{1}{1 - 2^{-\frac{1}{k-1}}} < \frac{dn }{16}.
$$
In the last inequality we used the facts that $2^{-x} \le 1-x/2$ for $0\le x \le 1$ and that $0<\frac{1}{k-1}\le 1$ for $k\ge 2$.
This is a contradiction as $e(G)= dn/2$. This proves the lemma.
\end{proof}

The Cauchy--Schwarz inequality together with 
a result by Hoory \cite[Lemma 1]{Hoory} yields the following lemma.

\begin{lemma}\label{lem: C2k}
Let $G$ be a bipartite graph with vertex partition $(A,B)$.
Suppose that the average of degrees of the vertices in $A$ is 
$d_A$
and the average of degrees of the vertices in $B$ is 
$d_B$. Then for every $k\in \mathbb{N}$ we have 
$$|\Hom(C_{2k},G)| \geq d_A^k \cdot d_B^k.$$
\end{lemma}

\section{Proofs of the main results}
Now we are ready to prove our 
main results.

\begin{proof}[Proof of Theorem~\ref{thm:even}]
By taking a subgraph, assume that $G$ is a bipartite graph with at least $50000 k^3 n^{1+1/k}$ edges and let $d \ge 10^5 k^3 n^{1/k}$ be the average degree of $G$. We further assume that $G$ has no subgraph with larger average degree, as otherwise we can just take that subgraph to be our graph.
Also assume that $G$ has no rainbow $2k$-cycle.
We apply Lemma~\ref{lem: regularization 1} to obtain a graph $G'$ and some $i\in \mathbb{N}$ such that 
\begin{itemize}
    \item $|A| = \frac{1}{k}2^{-\frac{ki}{k-1}}n$;
    \item $d_{G'}(a) \in [2^{i-6} d, 2^{i-5} d ]$ for all $a\in A$;
    \item $d_{G'}(b) \leq 4 d$ for all $b\in B$.
\end{itemize}
Indeed, we can obtain the equality in the 
first bullet point by deleting some vertices if necessary. 
Note that the first two conditions ensures that the average of degrees of vertices in $B$ is at least
$$\frac{2^{i-6}d|A|}{n} \geq \frac{d}{64k} 2^{-\frac{i}{k-1}}.$$
Apply Lemma~\ref{lem: C2k} to obtain that 
\begin{align}\label{eq: lower}
\Hom(C_{2k},G') \geq 2^{k(i-6)}d^k (\frac{d}{64k})^k 2^{-\frac{ki}{k-1}}.
\end{align}

As $d_{G'}(b)\leq 4d$ for all $b\in B$ and $\sum_{b\in B}d_{G'}(b)= e(G')$, the convexity of the function $f(x)=x^k$ yields that 
$$\sum_{b\in B} d_{G'}(b)^k \leq \frac{e(G')}{4d} \cdot (4d)^k \leq (4d)^{k-1}e(G').$$
Hence, Lemma~\ref{lem: degenerate} implies that 
\begin{align*}
    |\Hom(C_{2k},G')| &\leq (2k)^{2k} |\Hom(S_k,G')| \leq (2k)^{2k} \left(\sum_{a\in A} d_{G'}(a)^k+  \sum_{b\in B} d_{G'}(b)^k  \right) \\ &
    \leq (2k)^{2k}( |A| 2^{k(i-5)} d^k + (4d)^{k-1} e(G') ) \leq (2k)^{2k} (2^{k(i-5)} d^{k} + (4d)^{k-1}\cdot 2^{i-5}d)|A| \\& \leq (2k)^{2k} d^k (2^{k(i-5)}  + 4^{k-1}2^{i-5}) \cdot \frac{1}{k} 2^{-\frac{ki}{k-1}}n.
\end{align*}
Here, the penultimate inequality holds as 
$e(G')\leq 2^{i-5}d|A|$.
Combining this with \eqref{eq: lower}, we obtain  
$$ d^k < (10^5 k^3)^k n.$$
However, we assume that $d\ge 10^5 k^3 n^{1/k}$, a contradiction. This proves the theorem.
\end{proof}

\begin{proof}[Proof of Corollary~\ref{cor:supsat}] 
The proof proceeds as in Theorem~\ref{thm:even}. If at least half of the $\Hom(C_{2k},G)$ is non-degenerate, then we can bound it from below using~\eqref{eq: lower}; otherwise we reach a similar contradiction as now $|\Hom(C_{2k},G)|\le 2|\D(C_{2k},G)|$ and hence 
$|\Hom(C_{2k},G)| \leq 2\cdot (2k)^{2k} |\Hom(S_k,G)|$ instead of Lemma~\ref{lem: degenerate}.
\end{proof}

\begin{proof}[Proof of Theorem~\ref{thm:cycle}]
Suppose that $G$ does not have a rainbow cycle. By iteratively deleting low degree vertices, we may assume that the minimum degree of $G$ is $\delta \ge d(G)/2 \ge 32\log^2(5n)$. Apply Lemma~\ref{lem: regularization 2} to obtain a graph $G'$ on $n'\leq 4n$ vertices such that $G'$ doesn't contain a rainbow cycle and every vertex of $G'$ has degree between $\delta/2$ and $\delta$.

Let $k= \log n'$. Because $G'$ does not contain any rainbow $2k$-cycle, Lemma~\ref{lem: degenerate} implies that
$$ |\Hom(C_{2k},G')|\leq (2k)^{2k} |\Hom(S_k,G')|\leq (2k)^{2k} \delta^k n'.$$
On the other hand, even cycle satisfies Sidorenko's conjecture, so we know that 
$$(2k)^{2k} \delta^k n' \geq  (\delta/2)^{2k}.$$
As $n'=2^k$, this yields that 
$$ (4k^2\cdot \delta \cdot 2)^k \geq (\delta^2/4)^{k},$$ 
which is a contradiction as $\delta \ge 32\log^2(5n)>32k^2$. Hence $G$ must contain a rainbow cycle.
\end{proof}

\begin{proof}[Proof of Corollary~\ref{cor:3-graph}]
Consider a partition of $V(H)$ into $V_1,V_2,V_3$ where the number of edges containing one vertex from each is at least $\frac{1}{9} e(H)$. Consider an auxiliary bipartite graph $G$ with the vertex partition $(V_1, V_2)$ where $v_1v_2$ is an edge in $G$ with color $v_3$ if $v_1v_2v_3\in E(H)$ with $v_i\in V_i$. As $G$ contains at least $32n\log^2(5n)$ edges,  Theorem~\ref{thm:cycle} implies a rainbow cycle in $G$. This yields a loose cycle in $H$.
\end{proof}

It could be interesting to extend the above result to $3$-uniform hypergraphs that are close to linear, e.g. those with maximum co-degree being a constant, or even $o(n)$. One can also consider not necessarily proper colorings in which every vertex has bounded number of edges with the same color, or at most $o(n)$ edges of the same color.


\vspace{5mm}

\noindent\textbf{Acknowledgements.} While writing this note, Janzer and Sudakov \cite{JS22} obtained asymptotically the same bound $O(n \log^2 n)$ as our Theorem~\ref{thm:cycle} independently. Instead of regularization, they did a weighted homomorphism count. Their homomorphism inequalities are analogous to ours, and they are able to apply it to a much larger class of reflextive graphs, including e.g. hypercubes.

\bibliographystyle{abbrv}
\bibliography{reference}

\begin{thebibliography}{1}

\bibitem{BS74}
J.~A. Bondy and M.~Simonovits.
\newblock Cycles of even length in graphs.
\newblock {\em J. Combin. Theory Ser. B}, 16:97--105, 1974.

\bibitem{DLS13}
S.~Das, C.~Lee, and B.~Sudakov.
\newblock Rainbow {T}ur{\'a}n problem for even cycles.
\newblock {\em European J. Combin.}, 34(5):905--915, 2013.

\bibitem{ESi84}
P.~Erd\H{o}s and M.~Simonovits.
\newblock Cube-supersaturated graphs and related problems.
\newblock In {\em Progress in graph theory}, pages 203--218. Academic Press,
  Toronto, ON, 1984.

\bibitem{Hoory}
S.~Hoory.
\newblock The size of bipartite graphs with a given girth.
\newblock {\em J. Combin. Theory Ser. B}, 86(2):215--220, 2002.

\bibitem{J20}
O.~Janzer.
\newblock Rainbow {T}urán number of even cycles, repeated patterns and
  blow-ups of cycles.
\newblock To appear in Israel J. Math.

\bibitem{JS22}
O.~Janzer and B.~Sudakov.
\newblock On the {T}ur\'an number of the hypercube.
\newblock arXiv:2211.02015.

\bibitem{JS12}
T.~Jiang and R.~Seiver.
\newblock Tur\'an numbers of subdivided graphs.
\newblock {\em SIAM J. Discrete Math.}, 26:1238--1255, 2012.

\bibitem{KMSV07}
P.~Keevash, D.~Mubayi, B.~Sudakov, and J.~Verstra{\"e}te.
\newblock Rainbow {T}ur{\'a}n problems.
\newblock {\em Combin. Probab. Comput.}, 16(1):109--126, 2007.

\bibitem{T22}
I.~Tomon.
\newblock Robust (rainbow) subdivisions and simplicial cycles.
\newblock arXiv:2201.12309.

\end{thebibliography}
\end{document}